\documentclass{ifacconf}

\usepackage{graphicx}      
\usepackage{natbib}        
\usepackage{amsmath,amssymb}
\usepackage{stmaryrd}
\usepackage{tikz}
\usetikzlibrary{patterns}
\usetikzlibrary{decorations.pathreplacing}
\usepackage{pgf}
\usepackage{verbatim}
\usepackage{algorithm}
\usepackage[noend]{algpseudocode}

\newcommand{\R}{\mathbb{R}} 

\newcommand{\N}{\mathbb{N}}

\newcommand{\T}{\mathcal{T}}
\newcommand{\Div}[2]{\nabla _{#2} \cdot (#1)}

\newcommand{\overbar}[1]{\mkern 1.5mu\overline{\mkern-1.5mu#1\mkern-1.5mu}\mkern 1.5mu}

\newcommand{\Pacc}{P^{\text{ac}}_c ( \mathbb{R}^d \times \mathbb{R}^d )}\newcommand{\INTDom}[3]{\int_{#2} #1 \text{d} #3}
\newcommand{\normL}[3]{\parallel #1 \parallel _ {L^{#2}(#3)}}
\newcommand{\supp}[1]{\text{supp}(#1)}
\newtheorem{rmk}{Remark}
\newtheorem{Def}{Definition}

\begin{document}
\begin{frontmatter}

\title{Sparse control of kinetic cooperative systems to approximate alignment} 

\thanks[footnoteinfo]{This work has been carried out in the framework of Archim\`ede Labex (ANR-11-LABX-0033) and of the A*MIDEX project (ANR-11-IDEX-0001-02), funded by the ``Investissements d'Avenir" French Government programme managed by the French National Research Agency (ANR). \\ The authors acknowledge the support of the Grant ANR-16-CE33-0008-01 by ANR.}

\author[First]{Beno\^it Bonnet, Francesco Rossi} 

\address[First]{Aix Marseille Universit\'e, CNRS, ENSAM, Universit\'e de Toulon, LSIS, Marseille, France \\ {\tt benoit.bonnet@lsis.org} , {\tt francesco.rossi@lsis.org}}

\begin{abstract}                
Cooperative systems are systems in which the forces among agents are non-repulsive. The free evolution of such systems can tend to the formation of patterns, such as consensus or clustering, depending on the properties and intensity of the interaction forces between agents.

The \textit{kinetic} cooperative systems are obtained as the \textit{mean field limits} of these systems when the number of agents goes to infinity. These limit dynamics are described by transport partial differential equations involving non-local terms. 

In this article, we design a  simple and robust control strategy steering any kinetic cooperative system to approximate alignment. The computation of the control at each instant will only require knowledge of the size of the support of the crowd in the phase space and of the Lipschitz constant of the interaction forces.
Besides, the control we apply to our system is {\em sparse}, in the sense that it acts only on a small portion of the total population at each time. It also presents the features of being obtained through a \textit{constructive} procedure and to be independent on the number of agents, making it convenient for applications.
\end{abstract}


\end{frontmatter}

\section{Introduction}
\label{section:Intro}

The study of collective behaviour in systems of interacting agents has been the focus of a growing interest from several scientific communities during the past decades, e.g. in robotics (coordination of robots or drone swarms), in biology (crowds of animals), in sociology (information formation process) or in civil engineering (dynamical study of crowds of pedestrians, see e.g. \cite{CPT1}). In particular, it is well known that simple interaction rules between agents can promote the formation of global patterns. This phenomenom is usually referred to as \textit{self-organization} (see e.g. \cite{Bio_selforganization}). 
\\ However, the emergence of such patterns may be conditional to a certain number of hypotheses. For instance, a crowd with weak interactions will admit initial configurations for which no global self-organization can hope to arise. It is then natural to study whether it may be possible for an external action (e.g. a regulator) to enforce the formation of a pattern even in an unfavorable situation. This is the problem of \textbf{control of crowds}, which we shall address here in the particular case of kinetic cooperative systems. 

We recall the mathematical definition of \textit{cooperative models} in the finite-dimensional case, a well-known family of models used to describe crowds of interacting agents (see, e.g., \cite{angeli,smith}). Let us consider a set of $N \in \N^*$ interacting agents. In our case, the agents are supposed to be all identical and the dynamics of the $i$-th agent is given by
\begin{equation} \label{eq:Finite_dim_sys}
\dot x_i = v_i \hspace{0.1cm} , \hspace{0.1cm}  \dot v_i = \frac{1}{N}\sum_{j=1}^N \psi\left(x_j-x_i,v_j-v_i\right) ,
\end{equation}
where $\psi : \R^d \times \R^d \rightarrow \R^d$ is supposed to be non-negatively collinear to the velocity, namely
\begin{equation} \label{eq:Hyp1_psi}
\forall (x,v) \in \R^d\times \R^d \hspace{0.1cm} , \hspace{0.1cm} \psi(x,v) = \xi(x,v)v \hspace{0.1cm} , \hspace{0.1cm} \xi(\cdot,\cdot) \geq 0.
\end{equation}
We also assume that 
\begin{equation} \label{eq:Hyp2_psi}
\psi(\cdot,\cdot) \text{~ is $L$-Lipschitz}.
\end{equation}

In this article, we are interested in designing a control strategy for \textit{kinetic} cooperative systems. These systems are obtained as the \textit{mean field limits} of the finite-dimensional cooperative systems of the form \eqref{eq:Finite_dim_sys} when the number $N$ of agents goes to infinity. In this formalism, the crowd is represented at each time $t$ by its density $(x,v) \mapsto \mu(t,x,v)$. The time evolution of $\mu$ is described by the following transport Partial Differential Equation (PDE) with non-local terms 
\begin{equation} \label{eq:Transport_pde}
\partial_t \mu + v \cdot \nabla_x \mu + \Div{\Psi[\mu](x,v) \mu}{v} = 0 ,
\end{equation}
where $\Psi[\mu]:(x,v) \mapsto \INTDom{\psi(y-x,w-v)}{\R^d \times \R^d}{\mu(y,w)}$. 
This limit dynamics leads naturally to the design of control strategies that are independent on the number $N$ of agents. In practice these strategies can be applied to approximately control a finite-dimensional system containing a large number of agents, with an error that can easily be estimated as a function of $N$.

In the kinetic approach, the density $\mu$ of the crowd is modeled as a \textit{probability measure}. We introduce in this scope the definition of the functional spaces $P_\text{c}(\R^n)$ and $P_\text{c}^\text{ac}(\R^n)$, that are the natural setting to study our control problem (see e.g. \cite{EvansGariepy}).

\begin{Def}
The space $P_{\text{c}}(\R^n)$ is the set of all  probability measures on $\R^n$ with compact support, endowed with the weak topology of measures. \\ 
The space $P_{\text{c}}^\text{ac}(\R^n)$ is the subset of $P_{\text{c}}(\R^n)$ of all probability measures absolutely continuous with respect to the Lebesgue measure, i.e. the set of probability measures $\mu \in P_{\text{c}}(\R^n)$ for which there exists a Lebesgue-integrable function $f$ such that $\mu = f\text{d}x$ where $\text{d}x$ is the Lebesgue measure on $\R^n$. The function $f$ is called the \textit{density} of $\mu$ with respect to the Lebesgue measure.
\end{Def}
The definition of solutions for equations of the form \eqref{eq:Transport_pde} is then stated in terms of time-dependent curves in the space of probability measures. 
\begin{Def}
A solution $t \mapsto \mu(t)$ of \eqref{eq:Transport_pde} with initial datum $\mu^0 \in P_\text{c}(\R^n)$ is a curve in the space $P_\text{c}(\R^n)$ continuous with respect to time, satisfying \eqref{eq:Transport_pde} in the weak sense and such that $\mu(0) = \mu^0$.
\end{Def}
A natural idea to control \eqref{eq:Finite_dim_sys} is to add a control term $u_i$ to the dynamics of the $v_i$'s for all $i$ (see \cite{CaponigroFPT1}). Yet, for the mean field limit, all the agents of the crowd are supposed to be identical, thus,  one cannot impose a control localized specifically on one or several of these agents. \\
We are then compelled to introduce a space-dependent control of the form $(t,x,v) \mapsto \chi_{\omega(t)}u(t,x,v)$ where $\chi_{\omega(t)}$ is the indicator function of the time-dependent \textbf{control set} $\omega(t) \subset \R^n$. The controlled version of \eqref{eq:Transport_pde} then writes
\begin{equation} \label{eq:Controled_pde}
\partial_t \mu + v \cdot \nabla_x \mu+ \Div{[\Psi[\mu](x,v) + \chi_{\omega(t)}u(t,x,v)] \mu}{v} = 0 ,
\end{equation}
where $\chi_{\omega(t)}u(t,\cdot,\cdot)$ defines a Lipschitz vector field at all times $t \geq 0$. Furthermore, we impose our control to be \textit{sparse} and \textit{bounded}, i.e. that it can only act on a small portion of the crowd with limited amplitude at all times. These constraints write

$\diamond$ \textbf{Population sparsity constraint}
\begin{equation} \label{eq:Pop_const}
\mbox{$\INTDom{}{\omega(t)}{\mu(t)(x,v)} \leq c \hspace{1cm} \forall t \geq 0 ; $}
\end{equation}
$\diamond$ \textbf{Boundedness constraint} 
\begin{equation} \label{eq:Amp_const}
\normL{u(t,\cdot,\cdot)}{\infty}{\R^d \times \R^d} \leq 1 \hspace{0.4cm} \forall t \geq 0.
\end{equation}
%
%
In the sequel, we will be interested in the notion of \textit{approximate alignment}, which is defined as follows.
\begin{Def} \label{def:approx_alignment}
A solution $\mu \in C([0,+\infty),P_{\text{c}}(\R^d\times\R^d))$ to \eqref{eq:Transport_pde} or \eqref{eq:Controled_pde} is said to be $\epsilon$-approximately aligned around $v^* \in \R^d$ starting from time $T$ if $\supp{\mu(t)} \subset \R^d \times B(v^*,\epsilon)$ for any $t \geq T$.
\end{Def}
Our goal in this framework is to prove the following result. 
\begin{thm} \label{thm:eps_alignment}
Let $\mu^0 \in \Pacc$ be a given initial data for \eqref{eq:Controled_pde}. For any constant $c >0$, limit velocity $v^* \in \R^d$ and precision $\epsilon >0$, there exists a time $T$ and a Lipschitz-in-space control $u$ with support $\omega$ satisfying the constraints \eqref{eq:Pop_const} and \eqref{eq:Amp_const} such that the corresponding solution $\mu \in C([0,+\infty),P^{\text{ac}}_{\text{c}}(\R^d\times\R^d))$ of \eqref{eq:Controled_pde} is approximately aligned around $v^*$ with precision $\epsilon$ starting from time $T$.
\end{thm}

The structure of the paper is the following. We present in Section \ref{section:properties_kinetic} general notions concerning transport PDEs with non-local terms. We then prove Theorem \ref{thm:eps_alignment} in Section \ref{section:proof_thm} as follows: we introduce in Section \ref{subsection:Fundamental_step} the fundamental step of our control strategy and show in Section \ref{subsection:Procedure} how its iteration steers the dynamics \eqref{eq:Controled_pde} to approximate alignment.

\section{Kinetic cooperative systems}
\label{section:properties_kinetic}

\subsection{Transport PDEs with non-local velocities}

In this section we briefly introduce some notions and results concerning transport PDEs with non-local interactions of the form \eqref{eq:Transport_pde} and \eqref{eq:Controled_pde}. We first recall the definition of \textit{pushforward} of a measure by a Borel map and \textit{Wasserstein distance} (see more details in \cite{villani2}).
\begin{Def}
Given a Borel map $f : \R^n \rightarrow \R^n$, the pushforward of a probability measure $\mu$ defined on $\R^n$ through $f$ is the measure $f\# \mu$ satisfying:
\begin{equation} \label{eq:Pushforward}
f \# \mu(B) = \mu(f^{-1}(B)) 
\end{equation}
for any measurable subset $B \subset \R^n$.
\end{Def}
\begin{Def}
A transference plan $\pi$ between two probability measures $\mu, \nu \in P_\text{c}(\R^n)$ is a probability measure in $P_{\text{c}}(\R^{2n})$ which first and second marginals are respectively $\mu$ and $\nu$, namely, $\forall f,g \in C^{\infty}_{\text{c}}(\R^n), \INTDom{[f(x)+g(y)]}{\R^{2n}}{\pi(x,y)} = \INTDom{f(x)}{\R^n}{\mu(x)} + \INTDom{g(y)}{\R^n}{\nu(y)}$. We denote by $\Pi(\mu,\nu)$ the set of all transference plans between $\mu$ and $\nu$. The Wasserstein distance of order $p \geq 1 $ between $\mu$ and $\nu$ is then defined by 
\begin{equation*}
\mbox{$W_p(\mu,\nu) = \underset{\pi \in \Pi(\mu,\nu)}{\text{inf}} \left\{ \left( \INTDom{\vert x-y \vert ^p}{\R^{2n}}{\pi(x,y)} \right)^{\frac{1}{p}} \right\}.$}
\end{equation*}
\end{Def}
Due to its high convenience for computations and its numerous properties, the Wasserstein distance is a canonical object to study dynamics of probability measures. The fundamental result of existence and uniqueness for the general transport PDE with non-local terms
\begin{equation} \label{eq:General_pde}
\partial_t \mu + \Div{\Phi[\mu,t]\mu}{} = 0 ,
\end{equation}
is stated in terms of the Wasserstein distance (see \cite{ambrosio},\cite{Pedestrian}).
\begin{thm} \label{thm:PDE}
Assume that
$\Phi : \R \times P_\text{c}(\R^n) \rightarrow C^1(\R^n,\R^n)\cap L^{\infty}(\R^n,\R^n)$ satisfies the following properties:

$\bullet$ $\Phi[\mu,t](\cdot)$ is uniformly Lipschitz and with sublinear growth, i.e. there exist $L$,$M$ not depending on $\mu$ and $t$ such that $\vert \Phi[\mu,t](y) - \Phi[\mu,t](x) \vert \leq L\vert y - x \vert$ and $\vert \Phi[\mu,t](x)\vert \leq M(1+\vert x \vert)$ for any $x,y \in \R^n$. \\
$\bullet$ $\Phi[\mu,t]$ is a Lipschitz function with respect to $\mu$, i.e. there exists $K$ such that $\parallel \Phi[\mu,t] - \Phi[\nu,t] \parallel_{C^0} \leq KW_p(\mu,\nu)$ for any $\mu,\nu \in \times P^{\text{ac}}_\text{c}(\R^n)$. \\
$\bullet$ $\Phi[\mu,t]$ is measurable with respect to $t$. 

Then for any $\mu^0 \in P_\text{c}(\R^n)$, there exists a unique solution $\mu(\cdot) \in C([0,+\infty),P_\text{c}(\R^n))$ of \eqref{eq:General_pde}. Furthermore, the solutions of \eqref{eq:General_pde} depend continuously on their initial datum. \\
Let $(\T(t,\cdot))_{t \geq 0}$ be the flow of diffeomorphisms of $\R^n$ generated by the time-dependent vector field $\Phi[\mu,t](\cdot)$, defined as the unique solution of the Cauchy problem $\partial_t \T(t,x) = \Phi[\mu(t),t](\T(t,x)) ~ , ~ \T(0,x) = x$. Then, the solution $\mu(\cdot)$ of \eqref{eq:General_pde} with initial datum $\mu^0$ writes $\mu(t) = \T(t,\cdot)\#\mu^0$ for any $t \geq 0$. \\
In particular, $\mu^0 \in P_\text{c}^\text{ac}(\R^n)$ implies that $\mu(t) \in P_\text{c}^\text{ac}(\R^n)$ for all times $t \geq 0$.
\end{thm}
\begin{rmk} \label{rmk:Flow_diffeo}
Theorem \ref{thm:PDE} implies that $(\T(t,\cdot))_{t \geq0}$ describes the evoultion of the support of the measure $\supp{\mu(\cdot)}$. Indeed, for all times $t \geq 0$, any point $(x_t,v_t) \in \supp{\mu(t)}$ is the image of a corresponding point $(x_0,v_0) \in \supp{\mu^0}$ by the diffeomorphism $\T(t,\cdot)$.
\end{rmk}
In \eqref{eq:Transport_pde}, the vector field $\Phi[\mu,t]$ is $(x,v) \mapsto (v,\Psi[\mu])^T$. It can be easily checked that this vector field satisfies the hypotheses of Theorem \ref{thm:PDE}, see e.g. \cite{HaLiu1}. For this reason, we will define $\omega(\cdot)$ and $u(\cdot,\cdot,\cdot)$ such that our control $\chi_{\omega}u$ in \eqref{eq:Controled_pde} defines a Lipschitz vector field at all times, which ensures that the vector field $(x,v) \mapsto (v,\Psi[\mu]+\chi_{\omega}u)^T$ satisfies the hypotheses of Theorem \ref{thm:PDE}. 

We end this section by the statement of an estimate on the time evolution of the $L^{\infty}$-norm of the density of a probability measure following the dynamics \eqref{eq:Controled_pde}.
\begin{prop} \label{prop:Lp_estimates}
Let $\mu \in C([0,+\infty),\Pacc)$ be a solution of \eqref{eq:Controled_pde} with initial datum $\mu^0 \in \Pacc$ and $f(\cdot)$ be its density with respect to the Lebesgue measure. Then there holds for any times $t \geq 0$:
\begin{equation} \label{eq:Lp_estimates1} 
\begin{aligned}
& \frac{\text{d}}{\text{dt}} \parallel f(t) \parallel_{L^{\infty}} \leq \parallel f(t)\parallel_{L^{\infty}} \times \left[ \right. \parallel \Div{u(t)\chi_{\omega(t)}}{v} \parallel_{L^{\infty}} \\ 
& \hspace{2.05cm} + \parallel \Div{\Psi[\mu(t)]}{v} \parallel_{L^{\infty}(\supp{f(t)}} \left. \right] .
\end{aligned}
\end{equation}
\end{prop}
\textit{Proof:}
See \cite[Section 4.2]{controlKCS}. \hfill $\square$

\subsection{Invariance properties of kinetic cooperative systems}

We recall in this section the invariance properties of kinetic cooperative systems. One of the fundamental properties of \eqref{eq:Finite_dim_sys} is its invariance with respect to translations. Such properties are inherited by \eqref{eq:Transport_pde} and are stated as follows. 
\begin{prop} \label{prop:Invariance_trans}
Let $\mu(\cdot)$ be a solution of \eqref{eq:Transport_pde} with initial datum $\mu^0$, and $(y,w) \in \R^d \times \R^d$ a vector representing a translation. Define the curve $ \tilde{\mu}(t,x,v) = \mu(t,x+y+tw,v+w)$. Then $\tilde{\mu}(\cdot)$ is the unique solution of \eqref{eq:Transport_pde} with initial datum $\tilde{\mu}^0$, image of $\mu^0$ by the translation along $(y,w)^T$.
\end{prop}

Moreover, the attractivity of the interaction forces of \eqref{eq:Finite_dim_sys} allows us to establish an easy estimate of the evolution through time of the support of a solution of \eqref{eq:Transport_pde}.
\begin{prop} \label{prop:Invariance_supp}
Let $\mu(\cdot)$ be a solution of \eqref{eq:Transport_pde} with initial datum $\mu^0 \in \Pacc$. Then one has the following support invariance property: if $\supp{\mu^0} \subset \prod_{i=1}^d([\underbar{X}^i,\overbar{X}^i] \times [\underbar{V}^i,\overbar{V}^i])$ then $\supp{\mu(t)} \subset \prod_{i=1}^d([\underbar{X}^i+t\underbar{V}^i,\overbar{X}^i+t\overbar{V}^i] \times [\underbar{V}^i,\overbar{V}^i])$ for any $t \geq 0$.
\end{prop}

\textit{Proof:} This invariance is a direct consequence of Remark \ref{rmk:Flow_diffeo} and of the fact that $(x,v) \mapsto (v,\Psi[\mu](x,v))^T$ always points inward $\R \times [\underbar{V}^i,\overbar{V}^i]$ along $v$, for each $i \in \left\{1,..,N \right\}$. \hfill $\square$ 

We assume from now on that $\supp{\mu^0}$ is contained within the box $\prod_{i=1}^{d}([\underbar{X}^i,\overbar{X}^i]\times[\underbar{V}^i,\overbar{V}^i])$. The invariance properties given in Proposition \ref{prop:Invariance_trans} allow us to restrict the proof of Theorem \ref{thm:eps_alignment} to the case where $v^* = 0$, $\underbar{X}^i = 0$, $\overbar{X}^i >0$ and $\underbar{V}^i =0$, $\overbar{V}^i >0$ for any $i$, without loss of generality. Indeed, one can always achieve approximate alignment in the sense of Definition \ref{def:approx_alignment} in dimension $j$ with $\underbar{V}^j < 0 < \overbar{V}^j$ by applying the following strategy.

$(1)$ Define $\tilde{\mu}_j(t,x,v) = \mu(t,x-t\underbar{V}^je_j,v-\underbar{V}^je_j)$ (where $e_j$ stands for the $j$-th unitary vector of $\R^d$) and notice that it follows the dynamics \eqref{eq:Transport_pde} by Proposition \ref{prop:Invariance_trans}. Perform approximate alignment around $0$ with precison $-\underbar{V}^j + \frac{\epsilon}{\sqrt{d}}$. The velocity support of the orignal system in dimension $j$ is now $[\underbar{V}^j,\frac{\epsilon}{\sqrt{d}})$. \\
$(2)$ Define $\hat{\mu}_j(t,x,v) = \mu(t,-x-t\frac{\epsilon}{\sqrt{d}} e_j,-v-\frac{\epsilon}{\sqrt{d}}e_j)$ and notice that it satisfies the dynamics \eqref{eq:Transport_pde} with $\psi'(x,v) = - \psi(-x,-v)$. Perform approximate alignment around $0$ with precision $\frac{2\epsilon}{\sqrt{d}}$. The velocity support for the initial measure is included in $(-\frac{\epsilon}{\sqrt{d}},\frac{\epsilon}{\sqrt{d}})$ and approximate alignment is achevied for the $j$-th component.


\section{Proof of Theorem \ref{thm:eps_alignment}}
\label{section:proof_thm}

In this section we prove our main result Theorem \ref{thm:eps_alignment}, using a constructive algorithmic approach, in the spirit of \cite{controlKCS} and \cite{MTNS2016}. We assume henceforth that the dynamics \eqref{eq:Controled_pde} is unidimensional, i.e. $d = 1$. The case $d>1$ can be treated with a very similar technique, see e.g. \cite[Section 4.4]{controlKCS}. We define in Section \ref{subsection:Fundamental_step} the fundamental step of our control strategy, and then show in Section \ref{subsection:Procedure} how its iteration steers the dynamics to approximate alignment in the sense of Definition \ref{def:approx_alignment}.

\subsection{Fundamental step in 1D}
\label{subsection:Fundamental_step}

Assume that $\supp{\mu^0} \subset [0,X_0]\times[0,V_0]$. Our aim in this section is to build a Lipschitz control $\chi_{\omega(\cdot)}u^0(\cdot,\cdot,\cdot)$ and a time $T_0$ such that the control satisfies the constraints \eqref{eq:Pop_const} and \eqref{eq:Amp_const} at all times in $[0,T_0]$ and such that it reduces the size of $\supp{\mu(T_0)}$ along the velocity component. \\
%
To this end, we define a partition of our initial domain $[0,X_0]\times[0,V_0]$ into $ n = \lceil \frac{2}{c} \rceil$ rectangles $(\Omega_{[i]}^0)_{1 \leq i \leq n} = ([x_{[i-1]},x_{[i]}]\times[0,V_0])_{0 \leq i \leq n}$. The points $(x_{[i]})_{1 \leq i \leq n-1}$ are defined recursively to be the minimal values such that $\mu^0 \left( [x_{[i-1]},x_{[i]}] \times [0,V_0] \right) = \frac{c}{2}$ starting from $x_{[0]} = 0$. We also define $x_{[n]} = X_0$. \\
Note that $\mu^0 \left( [x_{[n-1]},x_{[n]}] \times [0,V_0] \right) \leq \frac{c}{2}$ and that the points $x_{[i]}$ are well defined, since $\mu^0 \in P^{\text{ac}}_{\text{c}}(\R\times\R)$, ensuring that $x \mapsto\mu^0([x_{[i-1]},x] \times [0,V_0])$ is a continuous function.
We further define the parameter $\delta_0$ as the biggest real number such that 
\begin{equation} \label{eq:delta_def}
\mu^0 ([x_{[i-1]}-3\delta_0,x_{[i]}+3\delta_0] \times \mathbb{R}) \leq c \text{~ $\forall i \in \left\{1,..,n \right\}$}.
\end{equation}
We fix a parameter $\eta_0 \in (0,V_0)$ and a time $T_0 \in \left( 0 , \frac{\delta_0}{V_0} \right]$, which precise choices shall be detailed in Section \ref{subsection:Procedure}. We define the control sets $(\omega_{[i]}^0)_{1 \leq i \leq n}$ by
\begin{equation}
\omega_{[i]}^0 = [x_{[i-1]}-2\delta_0 , x_{[i]}+2\delta_0] \times [0,V_0+\eta_0].
\end{equation}
We introduce the corresponding controls $(u_{[i]}(\cdot,\cdot,\cdot))_{1 \leq i \leq n}$ defined by $u_{[i]} : (t,x,v) \mapsto -\zeta_{[i]}(t,x,v)$ with $\zeta_{[i]}(\cdot,\cdot,\cdot)$ given by 
\begin{equation}
\zeta_{[i]}(t,x,v) = 
\left\{
\begin{aligned}
& 1 \hspace{0.2cm} \text{on} \hspace{0.2cm} [x_{[i-1]}-\delta_0 , x_{[i]}+\delta_0] \times [\eta_0,V_0] , \\ 
& \text{linearly decreasing to 0}  \hspace{0.2cm} \text{on} \hspace{0.2cm} \\&\hspace{0.2cm} \omega_{[i]}^0 \backslash ([x_{[i-1]}-\delta_0 , x_{[i]}+\delta_0] \times [\eta_0,V_0]), \\
& 0 \hspace{0.2cm} \text{on $(\R \times \R) \backslash \omega_{[i]}^0$}.
\end{aligned}
\right.
\end{equation}

A picture of this domain decomposition along with the definition of $\zeta_{[i]}(\cdot,\cdot,\cdot)$ for a given $i$ is given in Figure 1.

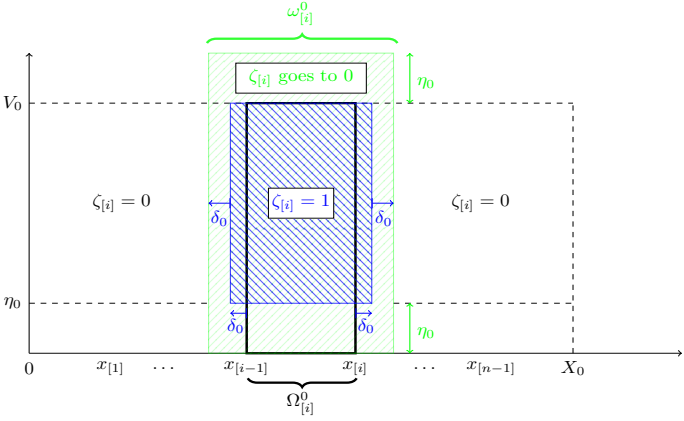
\begin{figure}[!h] \label{fig:Illustration}
\flushleft
\resizebox{0.5\textwidth}{0.23\textheight}{
\begin{tikzpicture}
\draw[->] (0,0)--(12,0);
\draw[->] (0,0)--(0,6); 
\draw(0,-0.3) node {0};
\draw(1.5,-0.3) node {$x_{[1]}$}; \draw(8.5,-0.3) node {$x_{[n-1]}$}; 
\draw(2.5,-0.3) node {$\dots$}; \draw(7.3,-0.3) node {$\dots$};
\draw[dashed] (0,5)--(10,5); \draw(-0.3,5) node [align = center] {$V_0$};
\draw[dashed] (10,0)--(10,5);\draw(10,-0.3) node [align = center] {$X_0$};
\draw[dashed] (0,1)--(3.3,1);\draw[dashed] (6.7,1)--(10,1);
\draw(-0.3,1) node [align = center] {$\eta_0$};
\draw[dashed] (4,1)--(4,0); \draw[dashed] (6,1)--(6,0);
\draw(4,-0.3) node {$x_{[i-1]}$}; \draw (6,-0.3) node {$x_{[i]}$};
\draw(5,-1) node {$\Omega^0_{[i]}$};
\draw[line width = 1.3pt](4,0) rectangle (6,5); 
\draw[blue, pattern=north west lines, pattern color=blue ,  opacity = 0.8] (3.7,1) rectangle (6.3,5);
\draw[green, pattern=north east lines, pattern color=green ,  opacity = 0.4] (3.3,0) rectangle (6.7,6);
\draw[fill = white] (4.4,2.7) rectangle (5.6,3.3); \draw[blue](5,3) node {$\zeta_{[i]} = 1$};
\draw[fill = white] (3.8,5.2) rectangle (6.2,5.8); \draw[green](5,5.5) node {$\zeta_{[i]}$ goes to 0};
\draw(1.7,3) node {$\zeta_{[i]} = 0$}; \draw(8.3,3) node {$\zeta_{[i]} = 0$}; 
\draw[line width = 1.3pt , decorate,decoration={brace,amplitude=5pt,mirror},xshift=2pt,yshift=0pt]
(3.95,-0.5) -- (5.95,-0.5);
\draw[line width = 1.3pt ,green!80, decorate,decoration={brace,amplitude=5pt},xshift=2pt,yshift=0pt]
(3.2,6.2) -- (6.7,6.2); \draw[green] (5,6.8) node {$\omega^0_{[i]}$};
\draw[|->,blue](6,0.8)--(6.3,0.8); \draw[blue](6.2,0.55) node {$\delta_0$};
\draw[|->,blue](4,0.8)--(3.7,0.8); \draw[blue](3.8,0.55) node {$\delta_0$};
\draw[|->,blue](6.3,3)--(6.7,3); \draw[blue](6.52,2.7) node {$\delta_0$};
\draw[|->,blue](3.7,3)--(3.3,3); \draw[blue](3.5,2.7) node {$\delta_0$};
\draw[<->,green](7,1)--(7,0); \draw[green](7.3,0.4) node {$\eta_0$};
\draw[<->,green](7,6)--(7,5); \draw[green](7.3,5.4) node {$\eta_0$};
\end{tikzpicture}
}
\caption{{\small \textit{Construction  of a control set $\omega_{[i]}^0$ and of the corresponding $\Omega_{[i]}^0$ and $\zeta_{[i]}(\cdot,\cdot,\cdot)$ for a given $i \in \left\{1,..,n \right\}$}}}
\end{figure}
We consider the time partition $[0,T_0] = \bigcup_{i=0}^{n-1} \left[ \frac{iT_0}{n}, \frac{(i+1)T_0}{n}\right]$ and apply each control $u_{[i]}(\cdot,\cdot,\cdot)$ on the set $\omega^0_{[i]}$ for $t \in \left[ \frac{iT_0}{n}, \frac{(i+1)T_0}{n}\right)$. This control design ensures the following properties.

$(1)$ The control is Lipschitz and satisfies \eqref{eq:Amp_const} at all times by definition of the functions $(\zeta_{[i]}(\cdot,\cdot,\cdot))_{1 \leq i \leq n}$. 

$(2)$  By Proposition \ref{prop:Invariance_supp}, one can easily check that $\supp{\mu(t)} \subset [0,X_0+tV_0]\times[0,V_0]$ for all times $t \in [0,T_0]$, yielding
\begin{equation} \label{eq:spatial_supp_evol}
\supp{\mu(t)} \subset [0,X_0+T_0V_0]\times[0,V_0]
\end{equation}
for all times $t \in [0,T_0]$. Moreover, choosing $T_0 \leq \delta_0/V_0$, we have by Remark \ref{rmk:Flow_diffeo} that all points in $\supp{\mu^0}$ will locally undergo a displacement of amplitude at most equal to $\delta_0$ in the variable $x$.

$(3)$ The population constraint \eqref{eq:Pop_const} is respected at all times. Indeed, for any $i \in \left\{1,..,n \right\}$ one has
\begin{equation}
\begin{aligned}
\INTDom{}{\omega_{[i]}^0}{\mu(t)(x,v)} \leq \int_{x_{[i-1]}-3\delta_0}^{x_{[i]}+3\delta_0} \int_0^{V_0} \text{d}\mu^0(x,v) = c.
\end{aligned}
\end{equation}

After having defined a proper control satisfying our constraints, we are interested in building estimates for the size of $\supp{\mu(T_0)}$. To do so, we monitor the evolution through time of the points $(x,v)$ such that $v$ realizes the maximum of velocity in $\Omega_{[i]}(t)$ (similar estimates were given in \cite{MTNS2016}). Here for $i \in \left\{1,..,n \right\}$ and $t\geq0$, we define $\Omega_{[i]}(t)$ as the image of $\Omega_{[i]}^0$ through the flow $\T_{u,\omega}(t,\cdot)$ generated by $(x,v) \mapsto (v,\Psi[\mu]+\chi_{\omega}u)^T$ as described in Theorem \ref{thm:PDE}. \\
We define the functions $b_i(t) = \text{sup} \left\{ v \hspace{0.1cm} \text{s.t.} \hspace{0.1cm} (x,v) \in \Omega_{[i]}(t) \right\}$. They satisfy the following properties.

$\diamond$ ~ If $b_i(t) \geq \eta_0 $ for all $t \in [0,T_0]$, then using \eqref{eq:Controled_pde}, Remark \ref{rmk:Flow_diffeo} and the Lipschitzianity of $\psi(\cdot,\cdot)$ one has that 

\begin{equation*}
\begin{aligned}
\dot b_i(t) & = \Psi[\mu(t)](x_i(t),b_i(t)) + \chi_{\omega_{[i]}^0}u_{[i]}(t,x_i(t),b_i(t)) \\
& \leq \mbox{$ L \int\limits_{\R\times\R} (w-b_i(t)) \text{d} \mu(t)(y,w) + \chi_{\omega_{[i]}^0}u_{[i]}(t,x_i(t),b_i(t)) $} \\
& \leq (V_0 - b_i(t)) + \chi_{\omega_{[i]}^0}u_{[i]}(t,x_i(t),b_i(t)).
\end{aligned}
\end{equation*}
Applying Gronwall lemma to $b_i(\cdot) - V_0$, noticing that $b_i(0)-V_0 \leq 0$ and taking $t = T_0$, we have
\begin{equation}
b_i(T_0) \leq V_0 + e^{-LT_0} \int_0^{T_0} \chi_{\omega_{[i]}^0}u_{[i]}(t,x_i(t),b_i(t)) \text{d}t.
\end{equation}
By construction of the set $\omega_{[i]}^0$ and of the control $u_{[i]}(\cdot,\cdot,\cdot)$, the fact that $b_i(\cdot) \geq \eta_0$ on $[0,T_0]$ implies that $u_{[i]}(\cdot,x_i(\cdot),b_i(\cdot))$ is equal to (-1) on $[\frac{iT_0}{n},\frac{(i+1)T_0}{n})$ and to 0 on $[0,\frac{(i-1)T_0}{n})\cup [\frac{(i+1)T_0}{n},T_0]$, leading to
\begin{equation}
b_i(T_0) \leq V_0 - \frac{e^{-LT_0}T_0}{n}.
\end{equation}  
\\
$\diamond$ ~If $b_i(t) < \eta_0 \hspace{0.1cm} \text{for some} \hspace{0.1cm} t \in [0,T_0] \hspace{0.1cm}$, define $\bar{t}$ to be the biggest time for which $b_i(\cdot) \leq \eta_0$. Notice then that $v-b_i(s) \leq V_0 - \eta_0$ for all $v \in [0,V_0]$ and $s \geq \bar{t}$. By a similar argument as in the previous point, one gets
\begin{equation}
\mbox{$b_i(T_0) \leq \eta_0 ( 1 - LT_0 ) + LV_0T_0.$}
\end{equation}
This holds in particular if $\bar{t} = T_0$, i.e. if $b_i(\cdot) \leq \eta_0$.

Since $(\Omega_{[i]}(t))_{1 \leq i \leq n}$ defines a covering of $\supp{\mu(t)}$ for all $t \geq 0$, these estimates together with \eqref{eq:spatial_supp_evol} yield $\supp{\mu(T_0)} \subset [0,X_1] \times [0,V_1]$ \text{with}
\begin{equation} \label{eq:support_final_est}
\left\{
\begin{aligned}
V_1 & = \text{max}\left\{ \right. V_0 - \frac{e^{-LT_0}T_0}{n} \hspace{0.1cm},\hspace{0.1cm} \eta_0 ( 1 - LT_0) + LT_0V_0 \left. \right\} , \\
X_1 & = X_0 + T_0 V_0 .
\end{aligned}
\right.
\end{equation}
\subsection{Proof of Theorem \ref{thm:eps_alignment} in 1D}
\label{subsection:Procedure}

In this section, we show how a sequence of fundamental steps as defined in Section \ref{subsection:Fundamental_step} (namely a sequence of choices of $\eta,T$) steers the system to approximate alignment. 

To this end, we will apply the following algorithm. 

\noindent\fbox{
    \parbox{0.91\linewidth}{
        \flushleft \textbf{Initialization :} \\
        \hspace{0.5cm} Let $\epsilon > 0$, $c > 0$ be given.
        \flushleft \textbf{Step k : Given the size $(X_k,V_k)$ of $\supp{\mu^k}$} \\
        	\hspace{1cm} \textbf{If $V_k < \epsilon$ :} \\
        	\hspace{1.2cm} Approximate alignment is already achieved. \\
        \hspace{1cm} \textbf{Else:} \\
        \hspace{1.2cm} Choose $T_k$ and $\eta_k$ as in \eqref{eq:Hyp_Tk_etak} \\
        \hspace{1cm} \textbf{End If} \\ 
         \flushleft \textbf{End :} Iterate while $V_k \geq \epsilon$ \\
    }
}

Fix $\epsilon >0$. We define for any $k \in \N$ the measure $\mu^{k+1} = \mu \left( \sum_{l=0}^{k}T_l \right)$. The estimate \eqref{eq:support_final_est} shows us that $\supp{\mu^{k+1}} \in [0,X_{k+1}]\times[0,V_{k+1}]$ where:
\begin{equation} \label{eq:Induction_def}
\left\{
\begin{aligned}
& V_{k+1} = \text{max}\left\{ V_k - \frac{e^{-LT_k}T_k}{n} \hspace{0.1cm} , \hspace{0.1cm} \eta_k ( 1 - LT_k) + LT_kV_k \right\}, \\
& X_{k+1} = X_k + V_k T_k .
\end{aligned}
\right.
\end{equation}
We build the corresponding partition $(\Omega_{[i]}^k)_{1 \leq i \leq n}$ of $[0,X_k]\times[0,V_k]$ and define the corresponding $\delta_k$ as in \eqref{eq:delta_def}. We also define the sets $(\omega^k_{[i]})_{1 \leq i \leq n}$, along with the corresponding controls $(u^k_{[i]})_{1 \leq i \leq n}$ as in Section \ref{subsection:Fundamental_step}. We set $\alpha = 1 +\frac{3}{nL\epsilon}$ and choose
\begin{equation} \label{eq:Hyp_Tk_etak}
\mbox{$T_k = \text{min} \left\{ \frac{\delta_k}{V_k} \hspace{0.1cm} , \hspace{0.1cm} \frac{1}{\alpha L} \right\} \hspace{0.1cm} , \hspace{0.1cm}\eta_k = \frac{1}{2} \left( V_k - \frac{e^{-LT_k}T_k}{n(1-LT_k)} \right)$}.
\end{equation}

We now want to show that the sequence $(V_k)_{k \in \mathbb{N}}$ defined above becomes smaller than $\epsilon$ \textbf{within a finite number of iterations} of our fundamental step. To do so, we  prove the slightly stronger result that $(V_k)_{k \in \mathbb{N}}$ converges to a limit $V_* < \epsilon/2$. This  implies that there exists $K \in \N$ such that $V_k < \epsilon$ for all $k \geq K$, hence that our algorithm stops. 

We first prove the following useful estimate.
\begin{lem} \label{lem:Lp_estimate}
Let $f^k$ be the density of $\mu^k$ with respect to the Lebesgue measure for a given $k \in \N$. Then, it holds
\begin{equation} \label{eq:Lp_estimates2}
12 \normL{f^k}{\infty}{\R^d\times\R^d} \delta_k V_k \geq c .
\end{equation}
where $V_k$ is the size of the support of $\mu^k$ along the velocity component and $\delta_k$ is defined as in \eqref{eq:delta_def}.
\end{lem}

\textit{Proof:} The proof follows from a simple geometric argument. Consider a given positive number $r$. Since the mass inside a set of the form $[x^k_{[i-1]},x^k_{[i]}]\times[0,V_k]$ is less or equal to $\frac{c}{2}$ for any $i \in \left\{1,..,n \right\}$ and $k \in \N$, then the mass contained in $[x^k_{[i-1]}-r,x^k_{[i]}+r]\times[0,V_k]$ is less than $\frac{c}{2} + 2rV_k \parallel f^k \parallel_{L^{\infty}}$. Besides, the mass of one of the sets $[x^k_{[i-1]}-3\delta_k,x^k_{[i]}+3\delta_k]\times[0,V_k]$ is equal to $c$, by definition of $\delta_k$. Taking $r = 3\delta_k$ yields the desired estimate. \hfill $\square$

We observe the following properties of our algorithm.

$(1)$ Choosing $T_k \leq \frac{1}{\alpha L}$ with $\alpha = 1 + \frac{3}{nL\epsilon}$ ensures that $\eta_k$ is strictly positive as long as $V_k \geq \epsilon/2$.

$(2)$ Choosing $T_k$ and $\eta_k$ as in \eqref{eq:Hyp_Tk_etak} ensures that
\begin{equation*}
\begin{aligned}
\eta_k = \frac{1}{2} \left( V_k - \frac{e^{-LT_k}T_k}{n(1-LT_k)} \right) & \leq V_k - \frac{e^{-LT_k}T_k}{n(1-LT_k)} \\
\Longrightarrow (1-LT_k)\eta_k + LT_kV_k & \leq V_k - \frac{e^{-LT_k}T_k}{n},
\end{aligned}
\end{equation*}
hence proving that $V_{k+1} = V_k - \frac{e^{-LT_k}T_k}{n}$ as long as $V_k \geq \epsilon/2$. This implies in particular that $(V_k)_{k \in \N}$ strictly decreases as long as  $V_k \geq \epsilon/2$.

$(3)$ Since $T_k \leq \frac{1}{\alpha L}$, then $V_k - V_{k+1} = \frac{e^{-LT_k}T_k}{n} \geq \frac{e^{-\frac{1}{\alpha}}T_k}{n}$ for any $k \in \N$ such that $V_k \geq \epsilon/2$. This implies that for any $K$ such that $V_K \geq \epsilon/2$ one has 
\begin{equation} \label{eq:Timeseries_bound}
V_0 \geq \frac{e^{-\frac{1}{\alpha}}}{n} \sum\limits_{k=0}^{K} T_k.
\end{equation}

%
\textbf{We now prove that our algorithm terminates in a finite number of iterations.
}

We prove it by contradiction. Assume that $V_k \geq \epsilon/2$ for all $k \in \N$. Then $(V_k)_{k \in \N}$ is strictly decreasing, bounded from below, and thus converges to a limit $V_* \geq \epsilon/2$. \\
This implies by \eqref{eq:Timeseries_bound} that $T_k \underset{k \rightarrow +\infty}{\longrightarrow} 0 $. Hence $\eta_k \underset{k \rightarrow +\infty}{\longrightarrow} V_*/2 > 0$ by \eqref{eq:Hyp_Tk_etak}. Thus we infer that there exists a constant $\bar{\eta}>0$ such that $\eta_k \geq \bar{\eta}$ for any $k \in \N$. 

By definition of $u^k(\cdot,\cdot,\cdot)$ one has for any $i \in \left\{1,..,n \right\}$ : 
\begin{equation*}
\normL{\Div{u_{[i]}^k(t)}{v}}{\infty}{\omega_{[i]}^0} 
\leq ~\frac{1}{\eta_k} ~ \leq ~\frac{1}{\bar{\eta}}.
\end{equation*}
Moreover, one has for any $(x,v) \in \mathbb{R}^d \times \mathbb{R}^d \hspace{0.1cm}:$
\begin{equation*}
\begin{aligned}
 & \vert \Div{\Psi[\mu(t)]}{v}(x,v)\vert  = \\
 & \vert -\sum_{k=1}^d \left[ \INTDom{\partial_{v_k} \psi(y-x,w-v)}{\mathbb{R}^d  \times \mathbb{R}^d}{\mu(t)(y,w)} \right] \vert \leq L d ,
\end{aligned}
\end{equation*}
by Lipschitzianity of $(x,v) \mapsto \psi(x,v)$. This leads to the following estimate for all $t \geq 0$:
\begin{equation} \label{eq:Inf_norm_estimate}
\begin{aligned}
\parallel \Div{u(t)\chi_{\omega(t)}}{v} \parallel_{L^{\infty}} + \parallel \Div{\Psi[\mu(t)]}{v} \parallel_{L^{\infty}} \leq \bar{F},
\end{aligned}
\end{equation}
where $\bar{F} = (Ld +1 / \bar{\eta} ) > 0$. 

Recall that $\bar{T} = \sum_{k=0}^{\infty} T_k$ is finite as a consequence of \eqref{eq:Timeseries_bound}. Then, combining \eqref{eq:Inf_norm_estimate} with \eqref{eq:Lp_estimates1} one has 
\begin{equation*}
\normL{f^k}{\infty}{\mathbb{R}^d \times \mathbb{R}^d} \leq \normL{f^0}{\infty}{\mathbb{R}^d \times \mathbb{R}^d} e^{\bar{F}\bar{T}} ~ \forall k \in \N.
\end{equation*}
This, together with \eqref{eq:Lp_estimates2} yields 
\begin{equation*}
\delta_k \geq \frac{c e^{-\bar{F}\bar{T}}}{12 \normL{f^0}{\infty}{\mathbb{R}^d \times \mathbb{R}^d} V_0} = \bar{\delta} > 0 \hspace{0.2cm} \forall k \in \mathbb{N}.
\end{equation*}
Finally, one has 
\begin{equation}
T_k = \text{min} \left\{ \frac{\delta_k}{V_k} \hspace{0.1cm} , \hspace{0.1cm} \frac{1}{\alpha L} \right\} \geq \text{min} \left\{ \frac{\bar{\delta}}{V_0} \hspace{0.1cm} , \hspace{0.1cm} \frac{1}{\alpha L} \right\} > 0
\end{equation}
for all $k \in \mathbb{N}$. This implies that  $\bar{T} = \sum_{k=0}^{+\infty}T_k$ diverges to infinity, since $(T_k)_{k \in \N}$ is uniformly bounded from below by a positive constant. This contradicts \eqref{eq:Timeseries_bound}. Hence, one concludes that $V_* < \epsilon/2$. Since the sequence $(V_k)_{k \in \N}$ is decreasing, we conclude that there exists $K \in \N$ such that $V_k < \epsilon$ for all $k \geq K$. Thus, the algorithm stops. We have proved that our control strategy steers the system to $\epsilon$-approximate alignment around 0 starting from time $\bar{T} = \sum_{k=0}^{K} T_k \leq e^{\frac{1}{\alpha}}\lceil \frac{2}{c} \rceil V_0$ for any given precision $\epsilon > 0$.

\bibliographystyle{ifacconf}
\bibliography{coopcontrol_IFAC}

\end{document}